\documentclass[a4paper,12pt]{amsart}
\usepackage{amsfonts, amsmath, amssymb, amscd, latexsym}

\input xy
\xyoption{all} \CompileMatrices

\newtheorem{thm}{Theorem}[section]
\newtheorem{prop}[thm]{Proposition}
\newtheorem{lem}[thm]{Lemma}
\newtheorem{cor}[thm]{Corollary}

\newtheorem{defn}[thm]{Definition}

\newcommand{\demo}{ {\it   Proof. }}

\begin{document}

%\title{A proof that all Seifert 3-manifold groups and all virtual surface groups
%are conjugacy separable}

\title[Conjugacy separability of Seifert Manifolds]{A proof that all Seifert 3-manifold
 groups and all virtual surface groups are conjugacy separable}

\author{Armando Martino}
\address{Centre de Recerca Matematica, Bellaterra, 08193, Spain}
\email{AMartino@crm.es}

%\subjclass[2000]{Primary ; Secondary }
%\thanks{}
%\keywords{}

\date{\today}

%\dedicatory{}

%%% ----------------------------------------------------------------------

\begin{abstract}
We prove that the fundamental group of any Seifert 3-manifold is conjugacy
separable. That is, conjugates may be distinguished in finite quotients or,
equivalently, conjugacy classes are closed in the pro-finite topology.
\end{abstract}

\maketitle

\section{Introduction}

There has been considerable interest in the separability properties
of groups, both from group theorists and topologists.

A subset $X$ of a group $G$ is called separable if for each $g
\not\in X$ there exists a map $\pi: G \to Q$, to a finite group $Q$
such that $g \pi \not  \in X\pi$.

It is well known that {\em residual finiteness}, where $X=\{ 1 \}$,
for a finitely presented group implies a solution to the word
problem. Similarly, {\em subgroup separability}, where $X$ is a
finitely generated subgroup, implies a solution to the membership
problem. And {\em conjugacy separability}, where $X$ is a conjugacy
class implies a solution to the conjugacy problem. Additionally,
these properties are related to the problem of lifting immersed
subspaces of topological spaces to embedded subspaces of finite
covers.

Recently, \cite{AKT} have proved that certain Seifert 3-manifold
groups are conjugacy separable and in this paper we prove that {\em all}
Seifert
3-manifold groups are conjugacy separable.

{\bf Theorem} {\em The fundamental group of any Seifert Fibered 3-manifold is conjugacy
separable.}

They are already known to be residually finite, \cite{hempel} and subgroup separable,
\cite{scott78}. Our approach is different from that in \cite{AKT}, in that we base our argument on
the algebraic structure of a Seifert 3-manifold group as an extension, rather than trying to
decompose it into amalgamated free products. One advantage of this approach is that it allows an
essentially unified treatment of these groups, although there is a certain distinction which
arises depending on whether the Seifert fibred space is `built' from a free group or the
fundamental group of a compact surface.

Our proof relies on the fact that virtually
surface groups are all conjugacy separable.
However, while this is already known that virtually free
groups are conjugacy separable, it was surprisingly unknown for
virtually surface groups. It is important to realise the contrast
between conjugacy separability and residual finiteness at this
point, since the latter is easily shown to pass to finite extensions
whereas for the former this is unclear and possibly untrue. However,
it is known that Fuchsian groups are conjugacy separable, and we
show how to argue from here to deduce that all virtually surface
groups are conjugacy separable.

\section{Background}

\subsection{Pro-finite topology}

Given a group, $G$, the {\em pro-finite topology} on $G$ has as a
basis all the cosets of finite index subgroups of $G$. Each such
coset is both open and closed in the pro-finite topology.

While the separability properties of $G$ can be described in terms
of the finite quotients of $G$, it is sometimes more convenient to
talk about the pro-finite topology instead.

For instance, $G$ is called {\em residually finite} if for any $1
\neq g \in G$, there exists a finite quotient of $G$ in which $g$
does not map to the identity element. In other words, there is a
normal subgroup $N$ of $G$ of finite index such that $g \not\in N$.
Equivalently, this means that $\{ 1 \}$ is closed in the pro-finite
topology of $G$. In fact, this is the same as saying that any one
element subset of $G$ is closed or, the seemingly stronger
statement, that $G$ is Hausdorff.

The group $G$ is called {\em subgroup separable} if for every
finitely generated subgroup, $H$ of $G$ and every $g \in G - H$
there exists a normal subgroup of finite index, $N$ of $G$ such that
$g \not\in HN$. Thus, there is a finite quotient of $G$ in which the
images of $g$ and $H$ are disjoint. As before, this is the same as
saying that every finitely generated subgroup is closed.

More generally, a subset $S \subseteq G$ is called separable if for
every $g \not\in S$, there is a finite quotient of $G$ in which $g$
and $S$ have disjoint images. Equivalently, $S$ is separable if it
is closed.

A group $G$ is called {\em conjugacy separable} if conjugacy classes
are separable.

If we have a subgroup $H$ of $G$, then there are two possible
topologies one can put on $H$. Namely, the subspace topology and the
pro-finite topology of $H$ itself. In general the subspace topology
may be more coarse, but not if $H$ has finite index.

\begin{lem}
Let $H$ be a finite index subgroup of $G$. Then the subspace
topology and the pro-finite topology of $H$ are the same.
\end{lem}
\demo If $K$ is a finite index subgroup of $G$ then $K x \cap H$ is
either empty, or is a coset of $K \cap H$, which has finite index in
$H$. Hence the subspace topology is contained in the pro-finite
topology.

Conversely, for any $K_0$ which is a finite index subgroup of $H$,
$K_0 x$ is open in $G$, since $K_0$ has finite index in $G$. Hence
the pro-finite topology is contained in the subspace topology.
$\Box$

\subsection{Seifert fibred spaces}

Throughout the paper we shall use the term {\em surface group} to
mean a group isomorphic to the fundamental group of a compact
surface. This includes the free abelian group of rank 2, the closed
hyperbolic surface groups as well as any finitely generated free
group. Since we are generally working up to finite index, we will
mainly restrict our attention to the orientable case. Using standard
terminology, we will call a group a {\em virtual surface group} if
it has a surface subgroup of finite index.

A Seifert fibred space is a 3-manifold that is `almost' a bundle.
This property can be made more precise when one considers the
fundamental group, $\pi_1(M)$ of a Seifert fibred space, $M$.

\begin{thm}
\label{ses} Let $M$ be a Seifert fibred space. Then there exists a
short exact sequence,
$$
1 \to C  \to \pi_1(M) \to H
$$
such that $C$ is cyclic and  $H$ contains a finite index subgroup isomorphic
to a surface group.
\end{thm}
\demo This is proved in \cite{scott83}, where it is shown that $H$ is the
fundamental group of a 2-orbifold. By \cite{scott83},
such a group is virtually a surface group. We note that it is sufficient
to consider the case where the 2-orbifold is without boundary,
by the remarks in \cite{scott83}, and since finitely generated subgroups of
surface groups are surface groups
(recall that we are including free groups as surface groups).
In fact, it turns out that if $\pi_1(M)$ is infinite -
 the only real case of interest for us -
then $C$ is also infinite. But we shall not need to use this fact. $\Box$

There is a rich literature dealing with the subject of Seifert
fibred spaces (\cite{JS}, \cite{seif}, for instance) and their
role in the theory of 3-manifolds. For our
purposes, however, the Theorem above provides sufficient information
on their fundamental groups in order to demonstrate conjugacy
separability.

\section{Virtual Surface Groups}

The goal of this section is to prove that any finitely generated
group having a surface subgroup of finite index is conjugacy
separable. Since virtually free and virtually free abelian groups
are already known to be conjugacy separable, we shall in fact only
be considering the case of closed hyperbolic surface groups. In
fact, since in \cite{FR} it has been proved that Fuchsian groups are
conjugacy separable, we shall show that this is sufficient for one
to deduce that {\em any} finite extension of a closed hyperbolic
surface group is conjugacy separable. The connection between the
general case and the Fuchsian case is provided by the Nielsen
Realisation Theorem.

First, however, we start with some terminology.

\begin{defn}
Given a group $S$ we call an automorphism, $\phi$, of $S$ {\em
virtually inner} if $\phi^n$ is inner for some integer $n$. Given a
virtually inner automorphism, $\phi$, let $S*_{\phi}$ denote the
group
$$
\langle S, t \ : t^n=x, t^{-1} g t = g \phi, \mbox{\rm for all } g
\in S \rangle,$$ where $n$ is chosen to be the least integer such
that $\phi^n$ is inner and $x$ is chosen so that  $g \phi^n= x^{-1}
g x$ for all $g \in S$. This group has $S$ as a subgroup of finite
index.

Also, given an automorphism $\phi$ of $S$, we consider the
equivalence relation, $\sim_{\phi}$ on $S$ where $g_1 \sim_{\phi}
g_2$ if and only if there exists an $h \in S$ such that $(h^{-1}
\phi) g_1 h = g_2$. We call this relation "twisted-$\phi$ conjugacy"
with "twisted-$\phi$ conjugacy classes". We denote the
twisted-$\phi$ conjugacy class of an element $h$ by $[h]_{\phi}$.
\end{defn}

The main results we shall use in the section are the following,

\begin{thm}[Nielsen Realisation Theorem, \cite{ker}, \cite{nielsen}]
Let $S$ be the fundamental group of a closed hyperbolic surface. Let
$H$ be a finite subgroup of $Out(S)$ and let $G$ be the pre-image of
$H$ in $Aut(S)$. Then $G$ has $Inn(S) \cong S$ as a subgroup of
finite index and is a Fuchsian group.
\end{thm}

\begin{thm}[\cite{FR}]
Fuchsian groups are conjugacy separable.
\end{thm}

Putting the above two together immediately produces the following
crucial corollary.

\begin{cor}
\label{finext} Let $S$ be the fundamental group of a closed
hyperbolic surface, and $\phi$ a virtually inner automorphism of
$S$. Then $S*_{\phi}$ is conjugacy separable.
\end{cor}
\demo Clearly, the image of $\phi$ in $Out(S)$ is a finite cyclic
subgroup whose pre-image in $Aut(S)$ is isomorphic to $S*_{\phi}$.
$\Box$

The purpose of introducing our seemingly cumbersome terminology
(twisted-$\phi$ conjugacy) is that it enables us to pass between
finite extensions of a group, $S$, by describing the separability
properties of the extension solely in terms of $S$.

\begin{prop}
\label{twisted} Let $S$ be a conjugacy separable group and $\phi$ a
virtually inner automorphism of $S$. If $S*_{\phi}$ is conjugacy
separable, then twisted-$\phi$ conjugacy classes in $S$ are closed
in the pro-finite topology of $S$.
\end{prop}
\demo Note that $S$ is a finite index subgroup of $S*_{\phi}$, so a
subset of $S$ is closed in the pro-finite topology of $S$ if and
only if it is closed in the pro-finite topology of $S*_{\phi}$.

Now consider the element, $tg$, for $g \in S$. Clearly, every
conjugate of $tg$ is a conjugate of $tg$ by some element of $S$.
Moreover, if $x \in S$ then,
$$
x^{-1} tg x = t (x^{-1} \phi) g x.
$$
Hence the conjugacy class of $tg$ in $S*_{\phi}$ is equal to
$$
t ([g]_{\phi}),
$$
where the twisted-$\phi$ conjugacy class is understood to be in $S$.
Now, since group multiplication is a homeomorphism, $[g]_{\phi}$
must be closed in $S*_{\phi}$ and hence in $S$. $\Box$

\begin{prop}
\label{cstwist} Let $G$ be a group which has $S$ as a normal
subgroup of finite index. Suppose that all twisted-$\phi$ conjugacy
classes of $S$ are closed in the pro-finite topology of $S$, for all
virtually inner automorphisms $\phi$. Then $G$ is conjugacy
separable.
\end{prop}
\demo Again, as $S$ has finite index, a subset of $S$ is closed in
$S$ if and only if it is closed in $G$.

Now consider a $g \in G$. Let $x_1, \ldots, x_r$ be a set of coset
representatives for $S$ in $G$ and let $g_i = x^{-1} g x_i$. Let
$\phi_i$ be the automorphism of $S$ induced by conjugation by the
element $g_i$. Clearly, each $\phi_i$ is virtually inner. We observe
the following,
$$
\begin{array}{rcl}
\{ w^{-1} g w \ : w \in G \} & = &
\bigcup_{i=1}^r \{ x^{-1} g_i x \ : \ x \in S \} \\
& = & \bigcup_{i=1}^r g_i [1]_{\phi_i}, \\
\end{array}
$$
where $[1]_{\phi_i}$ is understood to be a twisted-$\phi_i$
conjugacy class in $S$. By hypothesis, each $[1]_{\phi_i}$ is closed
in $S$ and thence in $G$. Thus each $g_i[1]_{\phi_i}$ is closed in
$G$ and thus so is the (finite) union. Thus we have shown that every
conjugacy class in $G$ is closed, which is another way of saying
that $G$ is conjugacy separable. $\Box$

It is now clear how to put these results together to get,

\begin{thm}
\label{virtcs}
 Let $G$ be a group having a surface subgroup of finite
index. Then $G$ is conjugacy separable.
\end{thm}
\demo Let $S$ be a surface subgroup of finite index which we assume,
without loss of generality, is normal in $G$. If $S$ is free this is
proved by \cite{dyer}. If $S$ is free abelian (the torus case), this
is proved by \cite{form}.

So we are left with the case where $S$ is the fundamental group of a
closed hyperbolic surface. Proposition~\ref{twisted} and
Corollary~\ref{finext} imply that all twisted-$\phi$ conjugacy
classes of $S$ are closed, for any virtually inner $\phi$.
Proposition~\ref{cstwist} then implies that $G$ is conjugacy
separable. $\Box$

\section{Extensions}

The main goal of this section is to prove Theorem~\ref{sep},
which says that finite-by-surface groups are surface-by-finite. That
is, they are virtually surface groups. We note that this is already
known, see \cite{kap}, exercise 4.7, but we felt that since it seems to be a
key observation which allows us to construct our proof of the fact
that Seifert 3-manifold groups are conjugacy separable, it would be
beneficial to include a proof of it. Additionally, the proof we present here
is, as far as we are aware, new and is elementary up to some well
known results concerning the residual properties of free groups and
nilpotent groups.

We start with the following lemma.
\begin{lem}
\label{order} Let $F$ be a finitely generated free group and $1 \neq g \in F$.
For any integer $n$ there exists a finite quotient, $Q_n$ of $F$, such that the
image of $g$ is a central element of order $n$ in $Q_n$.
\end{lem}
\demo It is clearly sufficient to prove the Lemma in the case where $n$ is a
power of some prime. So we consider an arbitrary prime $p$ and will show that
there is a quotient of $F$ in which the image of $g$ is central and has order
$p^k$ for any natural number $k$.

Now consider the central series for $F$. Write $\gamma_0=F$ and
$\gamma_{i+1}=[F, \gamma_i]$. As $F$ is residually nilpotent (\cite{neu})
we can find a
least $i$ such that $g \not\in \gamma_i$. Hence $g$ is non-trivial and central
in the free nilpotent group $F/\gamma_i$. The group $F/\gamma_i$ is a torsion
free nilpotent group, and hence is residually a finite $p$-group for all primes
$p$ (\cite{segal} or deduce this from \cite{neu} again). Thus we can find a
sequence of quotients $Q_k$ of $F$ such that each $Q_k$ is a $p$-group in
which the image of $g$ is central and $g$ has order $p^{n_k}$
for some unbounded sequence $n_k$. (Specifically, $Q_k$ is a finite $p$-group
quotient of $F/\gamma_i$ in which the images of $g^{p^i}$, for $1 \leq i \leq
k$, are non-trivial.)

However, we note that if $h$ is a central element of a finite $p$-group $Q$
with order $p^r$, then there is a quotient of $Q$ in which $h$ is central and
has order $p^{r-1}$. Namely, quotient out by all central elements of order $p$.
Hence we are done. $\Box$

\begin{lem}
\label{centralsplit} Consider the short exact sequence of groups,
$$
1 \to G \to \Gamma \to^{\pi} H \to 1
$$
where $G$ is a finite subgroup of the centre of $\Gamma$ and $H$ is a surface
group. Then $\Gamma$ has a finite index subgroup $\Gamma_0$ such that $\Gamma_0
\cap G$ is trivial.
\end{lem}
\demo This is trivial if $H$ is free so we assume that $H$ is a
one-relator group and we write $H=\langle X \ : \ r \rangle$. Now
let $F$ be the free group on the set $X$, and $\rho : F \to H$ the
natural epimorphism. By the universal property of free groups, there
is a map $\sigma : F \to \Gamma$ and a commuting triangle,
$$
\xymatrix{
& F \ar[dl]_{\sigma} \ar[dr]^{\rho} & \\
\Gamma \ar[rr]^{\pi} & &  H
}
$$
Now $r \sigma$ is clearly in the kernel of $\pi$ and hence is an element of
$G$. Suppose the order of $r \sigma$ is $n$. Then, by Lemma~\ref{order}, there
exists a finite quotient, $Q$, of $F$ in which the image of $r$ is central and
has order $n$. Let $F_0$ be the kernel of the map from $F$ to $Q$ and consider
the subgroup $\Gamma_0=F_0 \sigma$. Clearly, $\Gamma_0$ is a finite index
subgroup of $\Gamma$. We claim that $\Gamma_0$ is the required subgroup.

For consider an element $w \in \Gamma_0 \cap G$. Then $w = u  \sigma $ for some
$ u \in F_0$. Moreover, as $w$ lies in the kernel of $\pi$, $u$ must lie in the
kernel of $\rho$. Hence we can write $u$ as a product of conjugates of $r^{\pm
1}$ as follows,
$$
u = \prod_{i=1}^k g_i^{-1} r^{\epsilon_i} g_i$$ where each $\epsilon_i = \pm 1$
and the $g_i \in F$. Note that since $r$ is central in $F/F_0$ we have that
$u=r^m$ in $F/F_0$, where $m = \sum^k_{i=1} \epsilon_i$. In fact, since $r$ has
order exactly $n$ in $F/F_0$ and $u \in F_0$ we also have that $m$ is a
multiple of $n$.

However, since $r \sigma$ is central in $\Gamma$ we also get that $w=u \sigma=(r \sigma)^m$. As $r
\sigma$ has order $n$ and $m$ is a multiple of $n$, we deduce that $\Gamma_0 \cap G$ is trivial.
$\Box$

\begin{thm}
\label{sep} Consider a short exact sequence of groups,
$$
1 \to G \to \Gamma \to H \to 1,
$$
where $G$ is a finite group and $H$ is a surface group. Then
$\Gamma$ has a subgroup of finite index which is a surface group. In
particular, $\Gamma$ is conjugacy separable.
\end{thm}
\demo Note that $\Gamma$ acts on $G$ by conjugation. As this is a
finite group, the kernel of this action, $\Gamma_1$ is a finite
index subgroup of $\Gamma$ which centralises $G$. In particular, $G
\cap \Gamma_1$ is a finite central subgroup of $\Gamma_1$. Since the
image of $\Gamma_1$ in $H$ is a finite index subgroup of $H$ and
hence a surface group, we may apply Lemma~\ref{centralsplit} to
deduce that $\Gamma_1$ and hence $\Gamma$ has a subgroup of finite
index which intersects $G$ trivially. Since this subgroup is
isomorphic to a subgroup of finite index in $H$, we have shown that
$\Gamma$ is virtually a surface group. Hence, by
Theorem~\ref{virtcs}, $\Gamma$ is conjugacy separable. $\Box$

\section{Main argument}

Given a Seifert 3-manifold, $M$, we will show that $\pi_1(M)$ is
conjugacy separable by considering the short exact sequence given by
Theorem~\ref{ses}.
$$
1 \to C \to \pi_1(M) \to H \to 1,
$$
where $C$ is cyclic and  $H$ is virtually a surface group.
 In particular, we do this by
analysing the action on the normal subgroup $C$. In fact, by
Theorem~\ref{sep}, it is sufficient to consider the case when $C$ is infinite
cyclic, so henceforth we shall assume that $C \cong \mathbb{Z}$.

We let $h$
denote the generator of this $\mathbb{Z}$ and note that every
conjugate of $h$ is equal to either $h$ or $h^{-1}$. In particular,
the subgroup $\langle h^k \rangle$ is a normal subgroup of
$\pi_1(M)$ for each integer $k$. Our strategy will be to show that
if two elements in $\pi_1(M)$ are not conjugate, then they will also
fail to be conjugate in a quotient of $\pi_1(M)$ by some $\langle
h^k \rangle$. Since these quotients are all conjugacy separable by
the previous section, this is clearly sufficient.

The crux of the argument that follows is the content of the
following lemma.
\begin{lem}
\label{exp} Let $g \in \pi_1(M)$. Then there exist integers $\lambda,
\lambda_0$ so that
 $g$ conjugate to  $g h^n$ in $\pi_1(M)$ if and only if $n=k \lambda$ or $n=k\lambda+
 \lambda_0$,
 for some integer $k$.
\end{lem}
\demo Let $C$ be the pre-image in $\pi_1(M)$ of the centraliser of $g$ in
$\pi_1(X)$. That is, $C$ is precisely the subgroup of elements which conjugate
$g$ to an element of the form $g h^n$ for some integer $n$. Let $C^+$ be the
subgroup of $C$ of index at most $2$ which centralise $h$.

Consider $x_1, x_2 \in C^+$ and suppose that $x_i^{-1} g x_i = g h^{n_i}$. Then
it is easy to check that
$$
x_2^{\mp 1} x_1^{-1} g x_1 x_2^{\pm 1} = g h^{n_1 \pm n_2}.
$$
In other words, the set $\{ h^n \ : x^{-1}gx=g h^n \ \mbox{\rm for some } x \in
C^+ \}$ is a subgroup of $\langle h \rangle$. Thus it has a generator
$h^{\lambda}$ and we have shown that $g h^n = x^{-1} g x$ for some $x \in C^+$
if and only if $n = k \lambda$.

If $C=C^+$ we are done by setting $\lambda_0=0$. Otherwise we can find a $y \in
C - C^+$ and an integer $\lambda_0$ such that $y^{-1} g y = g h^{\lambda_0}$.
The conclusion of this lemma is now immediate on noting that,
$$
y^{-1} g h^{k \lambda} y = g h^{-k \lambda + \lambda_0},$$ and that $C=C^+ \cup
C^+ y$. $\Box$

\begin{thm}
All Seifert 3-manifold groups are conjugacy separable.
\end{thm}
\demo By Theorem~\ref{sep}, it is sufficient to prove the theorem in the case
where we have a short exact sequence of groups associated
to a Seifert 3-manifold, $M$, of the form
$$
1 \to \mathbb{Z} \to \pi_1(M) \to H \to 1,
$$
where $H$ is a virtual surface group and $\mathbb{Z}$ is generated
by an element $h$.

So let us consider two elements $g, g_1 \in \pi_1(M)$ which are not
conjugate. In order to show that these elements are not conjugate in
some finite quotient, it is clearly sufficient to show that they are
not conjugate in some conjugacy separable quotient of $\pi_1(M)$.
Thus we may assume that the images of $g$ and $g_1$ {\em are} conjugate in
$H$, since $H$ is a conjugacy separable quotient of $\pi_1(M)$, by
Theorem~\ref{virtcs}. After taking suitable conjugates, we may
further assume that $g$ and $g_1$ have the same image in $H$.

Thus $g_1=g h^n$ for some integer $n$. However, by Lemma~\ref{exp}, there exist
integers $\lambda, \lambda_0$ so that an element of the form $g h^m$ is
conjugate to $g$ in $\pi_1(M)$ if and only if $m=k \lambda$ or $m=k \lambda +
\lambda_0$ for some integer $k$. In particular, since $g$ and $g_1$ are not
conjugate in $\pi_1(M)$, neither can they be conjugate in the quotient
$\pi_1(M)/\langle h^{\lambda} \rangle$. However, again by
Theorem~\ref{sep}, $\pi_1(M)/\langle h^{\lambda} \rangle$ is conjugacy
separable and we are done. $\Box$

\section*{Acknowledgments}

We gratefully acknowledge the support of the CRM and its staff and thank everyone at the CRM for
their hospitality. We are also very grateful to Noel Brady for many helpful conversations while
writing this work, and we would like to thank all the visitors at the CRM for a stimulating
environment in 2004-2005.

%\bibliographystyle{amsplain}
%\bibliography{refs}

\providecommand{\bysame}{\leavevmode\hbox to3em{\hrulefill}\thinspace}
\providecommand{\MR}{\relax\ifhmode\unskip\space\fi MR }
% \MRhref is called by the amsart/book/proc definition of \MR.
\providecommand{\MRhref}[2]{%
  \href{http://www.ams.org/mathscinet-getitem?mr=#1}{#2}
}
\providecommand{\href}[2]{#2}

\end{document}